\pdfoutput = 1

\documentclass[10pt,letterpaper]{article}
\usepackage[top=0.85in,left=2.75in,footskip=0.75in]{geometry}

\usepackage{changepage}

\usepackage[utf8x]{inputenc}

\usepackage{textcomp,marvosym}

\usepackage{fixltx2e}

\usepackage{amsmath,amssymb}

\usepackage{cite}

\usepackage{nameref}

\usepackage[right]{lineno}

\usepackage{microtype}
\DisableLigatures[f]{encoding = *, family = * }

\raggedright
\usepackage[aboveskip=1pt,labelfont=bf,labelsep=period,justification=raggedright,singlelinecheck=off]{caption}

\makeatletter
\renewcommand{\@biblabel}[1]{\quad#1.}
\makeatother

\date{}

\usepackage{lastpage,fancyhdr,graphicx}
\usepackage{epstopdf}
\fancyheadoffset[L]{2.25in}
\fancyfootoffset[L]{2.25in}
\usepackage{xcolor}

\begin{document}
\vspace*{0.2in}

\begin{flushleft}
{\Large
\textbf\newline{Gender Representation on Journal Editorial Boards in the Mathematical Sciences} 
}
\newline
\\
Chad M. Topaz\textsuperscript{1,\Yinyang,*},
Shilad Sen\textsuperscript{1\Yinyang}
\\
\bigskip
\textbf{1} Dept. of Mathematics, Statistics, and Computer Science, Macalester College, Saint Paul, MN, USA
\\
\bigskip

%
%
\Yinyang These authors contributed equally to this work.

* ctopaz@macalester.edu

\end{flushleft}
\section*{Abstract}
We study gender representation on the editorial boards of 435 journals in the mathematical sciences. Women are known to comprise approximately 15\% of tenure-stream faculty positions in doctoral-granting mathematical sciences departments in the United States. Compared to this pool, the likely source of journal editorships, we find that 8.9\% of the 13067 editorships in our study are held by women. We describe group variations within the editorships by identifying specific journals, subfields, publishers, and countries that significantly exceed or fall short of this average. To enable our study, we develop a semi-automated method for inferring gender that has an estimated accuracy of 97.5\%. Our findings provide the first measure of gender distribution on editorial boards in the mathematical sciences, offer insights that suggest future studies in the mathematical sciences, and introduce new methods that enable large-scale studies of gender distribution in other fields.


\section*{Introduction}
\label{sec:introduction}

We study gender representation on the editorial boards of 435 mathematical sciences journals, comprising over 13000 editorships. Our study serves three purposes. First, it provides a snapshot of gender on these editorial boards, which to our knowledge, has not before been measured. Second, it provides a benchmark to which future measurements can be compared, thus enabling longitudinal assessments of any changes over time. And third, it presents a methodology for large-scale gender representation studies that is scalable and largely automatic, thus facilitating those future investigations.

\begin{figure}[!h]
\includegraphics[width=\textwidth]{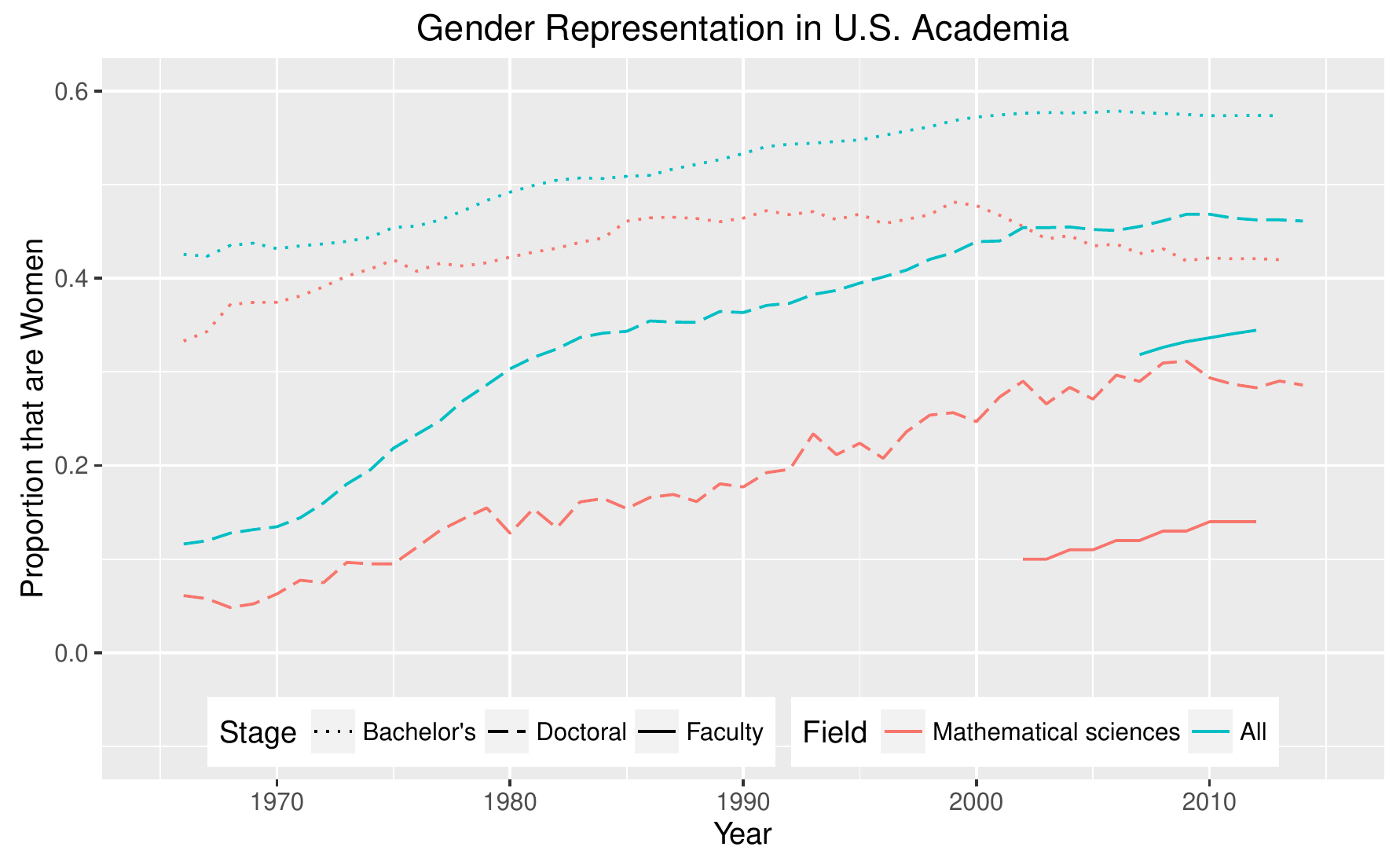}
\caption{{\bf Representation of women at various levels of academia in the United States.} Dotted blue: Bachelor's degree recipients in all fields \cite{IPEDSCompletions2016}. Dotted red: Bachelor's degree recipients in the mathematical sciences \cite{IPEDSCompletions2016}. Dashed blue: Doctoral degree recipients in all fields \cite{IPEDSCompletions2016,SED2016}. Dashed red: Doctoral degree recipients in the mathematical sciences \cite{IPEDSCompletions2016,SED2016}. Solid blue: Tenured/tenure-eligible faculty in all fields, restricted to doctoral granting institutions\cite{IPEDSSalaries2016}. Only limited data is available because \cite{IPEDSSalaries2016} did not gather tenure status data during many recent years. Solid red: Tenured/tenure-eligible faculty in the mathematical sciences, restricted to doctoral granting institutions \cite{AMS2016}.}
\label{fig1}
\end{figure}

From bachelor's degree students to tenured research faculty, women in higher education are grievously underrepresented within the mathematical sciences; see Fig~\ref{fig1}. In the United States, women comprise approximately 51\% of the population \cite{Census2016}. The percentage of bachelor's degrees in mathematics and statistics given to women has been falling since 1999 \cite{IPEDSCompletions2016}. In 2013, it reached 42\%, equal to that in 1979, revealing a setback of over more than three decades. At the doctoral level, despite gains for women in overall doctoral degrees granted, the proportion in the mathematical sciences has stagnated over the past decade, fluctuating around 29\% \cite{SED2016,IPEDSCompletions2016}. While the representation of women among faculty at doctoral granting institutions has been improving overall for the past decade, the proportion of mathematical sciences faculty positions held by women in 2013 was still merely 15\% for tenured/tenure-eligible faculty \cite{IPEDSSalaries2016, AMS2016}.

Beyond our understanding of the many historical barriers to women entering the workforce, there is a robust scholarly literature that examines the possible causes of additional underrepresentation in science, technology, engineering, and mathematics (STEM) disciplines. Evidence strongly suggests that males and females have similar inherent capacity in these fields \cite{Spe2005,HydLin2006}. With innate biological differences in ability excluded, one must ask what other factors contribute. We mention a few key studies. 

Work in \cite{LesCimMey2015} concludes that representation of women is poorer in fields whose practitioners believe in the importance of innate ``talent'' or ``brilliance'' to success in that field. In fact, in this study, mathematics comes in second place (behind philosophy) for field-specific ability belief, a measurement of how much emphasis that field's practitioners place on brilliance. A follow-up study finds that the terms ``brilliant'' and ``genius'' are attributed to men in STEM fields at much higher rates than they are to women. Furthermore, the frequency with which these terms appear in ratings of faculty on the RateMyProfessor.com teaching evaluations web site predicts, to some extent, diversity within a given academic field \cite{StoHorCim2016}. In another stark example of bias, scientists at research universities were asked to evaluate student application materials. Participants rated the same materials more highly when those materials carried a man's name than when they carried a woman's \cite{MosDovBre2012}. In the realm of prizes and awards in STEM fields, a gap between the proportions of women nominees and women prizewinners suggests an unconsciously or consciously biased selection process \cite{LinPinKos2012,PopLeb2012}. More comprehensively, a study solicited by the National Science Foundation and conducted by the American Association of University Women identifies many complex, interacting factors that contribute to underrepresentation of women, including: stereotypes about women's abilities; harsher self-assessment of scientific ability by women than by men; academic and professional climates that are dissatisfying to women; and unconscious bias \cite{Whysofew2010}.

In our present study, we examine gender representation on editorial boards of mathematical sciences journals. Given that women are underrepresented at so many levels in the mathematical sciences, one might ask why editorial boards merit special focus. Most straightforwardly, the representation of women on editorial boards could serve as an indicator for the field at large. However, we offer three additional reasons. First, editors are in positions of power. They act as stewards who ensure the quality of research being reported and whose editorial decisions even influence the broad course of research within the community; see \cite{ChoJohSch2014} and references therein. If women are underrepresented on editorial boards, it means that the field is being deprived of women's contributions and perspectives \cite{AddVil2003}. Relatedly, groups with more diversity are associated with making better decisions \cite{LoyWanPhi2013}, suggesting that a diverse journal editorial board benefits the scientific community that reads the journal. Second, editorial board positions provide valuable opportunities for intellectual growth and even for professional networking, which fosters career development and can contribute positively to tenure and promotion decisions \cite{GerDutBar2000,MetHar2009,ChoJohSch2014}. Women excluded from editorial board membership do not have access to these advantages. Finally, we note the importance of the role model effect, meaning that the presence of women on editorial boards might encourage more junior women to remain in academia \cite{SmiErb1986,BetLon2005,StePalAss2011}. In summary, a greater representation of women on editorial boards is to the benefit of the research community at large, to individual women, and to the future. 

Previous studies of gender representation on editorial boards have been carried out in medicine \cite{MorSon2007,AmrLanFah2011}, environmental biology and resource management \cite{ChoJohSch2014}, ecology \cite{FoxBurMey2016}, science (broadly) \cite{MauHilMor2013}, economics \cite{AddVil2003}, management \cite{MetHar2009}, and political science \cite{StePalAss2011}, to name a few examples.
These studies typically require a great deal of labor. Investigators must determine editorial board membership, and in an even more difficult task, they must infer the gender of each board member. To our knowledge, the aforementioned studies infer each editor's gender manually; the investigators use personal knowledge, make personal judgments of an editor's name, inquire with others, seek textual references to the editor using gender-specific pronouns, search for photos, and so forth. 

A few studies of gender representation in academia, with foci other than editorial boards, have taken automated approaches that infer gender based on first names. In the realm of article reviewers, \cite{FoxBurMey2016} infers the gender of a subset of its database of 8533 individuals by using a classification tool that compares each first name to a large database of known names. Similarly, in the realm of article authors, \cite{BriWhi2015} adopts a purely automated approach. In this innovative study, the investigators consider 100183 authors of mathematics papers posted to the arXiv scientific preprint repository over a period of 23 years and examine the association between gender representation and mathematical subfield. They leverage the previous work of \cite{Pie2014}, which inferred gender by comparing to a database of 40000 names classified by native speakers. While these studies reduce manual effort, they lack coverage for rare names, and can perform inaccurately for gender-ambiguous names. For example, in \cite{Pie2014}, gender could not be determined for over 40\% of the authors. 

Our own data-gathering and gender inference methods combine the strengths of the manual and automatic name-based approaches outlined above. However, instead of relying on manual work from the principal investigators, our study uses Amazon Mechanical Turk (MTurk), a crowdsourcing internet marketplace, and puts into place a variety of safe-guards to ensure high-quality, representative data. Advantages of our process include that (1) it requires minimal investigator intervention, (2) it exhibits high coverage and accuracy in gender inference, estimated at 99\% and 97.5\% respectively in our study, and (3) it reduces systemic bias that might erroneously inflate the number of reported men or women. With these benefits, our methodology is well-positioned to support future large-scale gender studies. For example, researchers could apply our method to conduct longitudinal analyses that require repeated studies of gender representation. Equally, one could apply our methods to larger fields and sample populations.

The rest of this paper is organized as follows. In the following section, we present our methods, including how the MTurk platform and automated gender inference tool work, how we gather and clean our data, and how we correct for bias in our gender inference procedure. We also provide a brief discussion of the complex construct that is gender, and some of the pitfalls of inferring another individual's gender. Then, we proceed to our results. Of 13067 editorships, 90.3\% are held by men, 8.9\% are held by women, and 0.8\% could not be determined. Of the 435 journals, 51 have no editorships held by women, and the median journal board includes just 7.6\% women. Additionally, we examine the associations between gender and journal subfield, impact factor, publisher, title on editorial board, and country. Finally, we conclude by reflecting on our results and offering some directions for further investigation.

\section*{Methods}
\label{sec:methods}

\begin{figure}[!h]
\includegraphics[width=\textwidth]{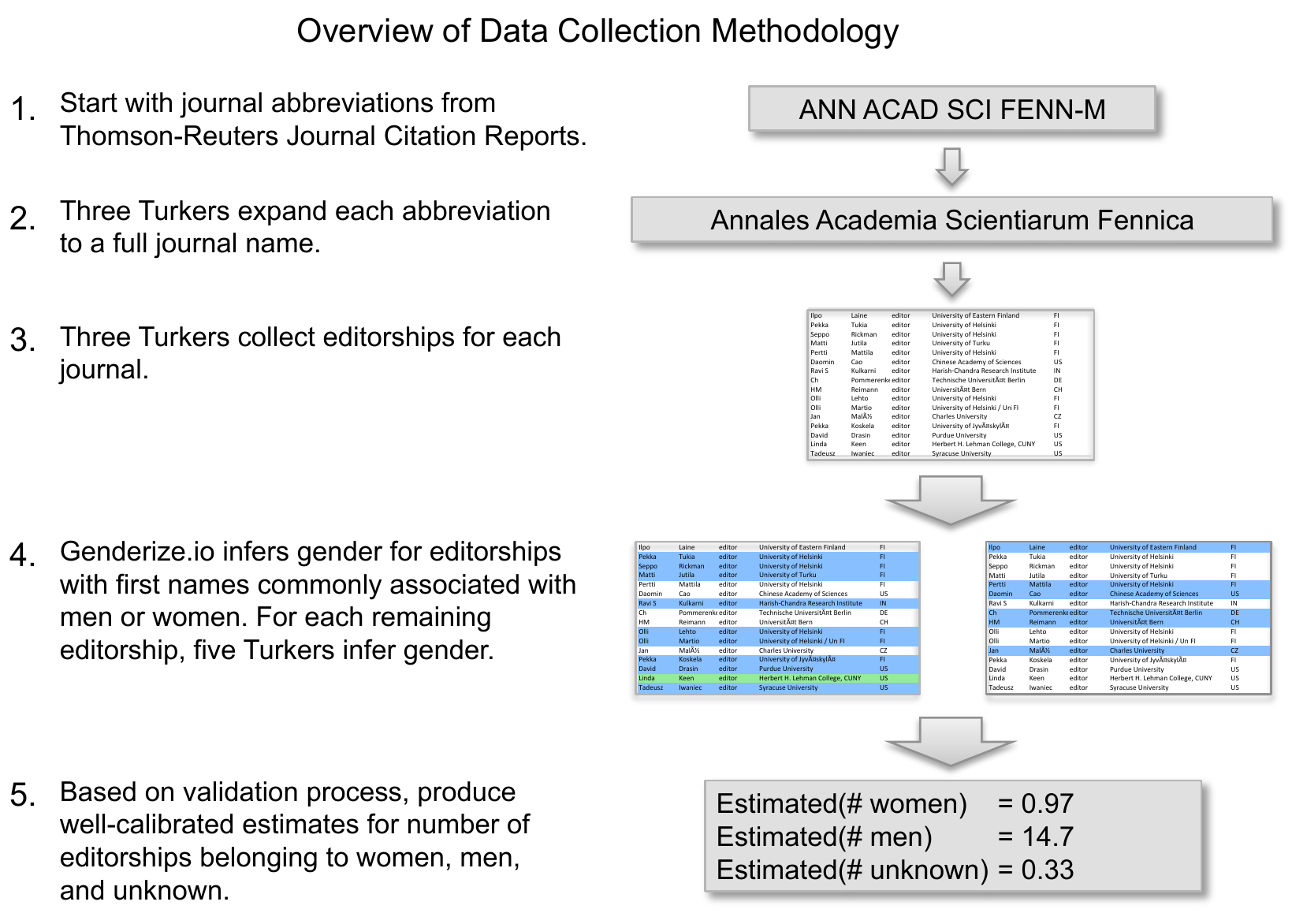}
\caption{{\bf Overview of data collection methodology used in this paper}. We combine publicly available data, algorithmic inference, crowdsourced data from Amazon's Mechanical Turk, and manual validation procedures. In the snapshots of data shown above, blue and green indicate editorships inferred to be held by men and women, respectively. For the example shown, the raw inferences are 15 men and one woman. A calibration procedure, described at length in the text, produces revised values of 14.7 men, 0.97 women, and 0.33 unknown.}
\label{fig2}
\end{figure}

Momentarily, we will explain our research methodology at length. First, to provide an overview, Fig~\ref{fig2} summarizes our procedures. We began with a collection of 605 abbreviations from journals in the mathematical sciences from an established source. Crowdworkers from Amazon Mechanical Turk then expanded those abbreviations into full journal names and collected the individual editorships associated with each journal. For editorships whose first names were overwhelmingly associated with being masculine or feminine, we used an automated tool to infer gender based solely on the name. For the remaining editorships, we asked the crowdworkers to use a search engine to infer gender. Finally we applied a validation and calibration procedure to estimate the accuracy of our approach and minimize structural prediction bias that would infer too many men or women.

\subsection*{Data collection}

There exist several indices of mathematical sciences journals, including Google Scholar, the American Mathematical Society's MathSciNet, and Thomson-Reuters Journal Citation Reports (JCR). While recognizing that any of these might serve as an appropriate starting point, we opted to use JCR because of its broad coverage and its impact factor data \cite{JCR2015}. From JCR, we obtained a list of 605 journals in the mathematical sciences. The JCR product we accessed provides only abbreviated journal titles, so a later step in our data collection involved determining the full titles of these journals. Within JCR, there are four categories of mathematical sciences journals: \emph{mathematics}, \emph{mathematics applied}, \emph{mathematical and computational biology}, and \emph{mathematics interdisciplinary applications}. Some journals appear in more than one category. For each journal, in addition to its abbreviated title, we saved its ISSN number (a code used to identify serial publications) and its five-year impact factor.

Journals store their editorial board information in a wide variety of online formats, ranging from highly-structured HTML web pages to .pdf documents and more. Thus, it is impractical to write computer code to gather this data. Equally, it is impractical for a small team of researchers to gather it manually for hundreds of journals. Instead, to locate online information about journal editorial boards, we used Amazon Mechanical Turk (MTurk). MTurk is an online marketplace for Human Intelligence Task (HIT) labor. HITs refer to tasks that cannot (or cannot easily) be performed by computers. Classic examples of HITs include: classifying an online product review as positive or negative, generating text descriptions of an image, and transcribing an audio clip. In the MTurk system, a \emph{requester} posts HITs which MTurk workers, known as \emph{Turkers}, may browse and choose to complete for a fee set by the requester. Once the Turker completes the task, the requester may review the work and approve or reject it. MTurk has seen growing use in the social sciences and in computer science \cite{BuhKwaGos2011,SenGieGol2015}.

Throughout our study, to ensure a high quality of data, we put in place four safeguards. First, we required that our Turkers had completed at least 1000 HITs in the past so that they would be sufficiently experienced. Second, we required that our Turkers had an approval rate of 99\% or higher for their past work. Third, we had each of our HITs completed by multiple Turkers in order to validate the data. Finally, we manually validated subsamples of the data.

For the 605 journals in our original data set, we asked Turkers to provide the full journal title as well as a link to online information about the journal's editorial board. Each HIT was completed by three Turkers. For 406 journals, all respondents reported the same full title. For the remaining 199 journals, we resolved the full title manually. For 183 editorial board URLs, all respondents reported the same URL. For the remaining 422 URLs, we resolved the URL manually. Two of the journals were clearly defunct. We removed these from our data set, resulting in 603 journal titles and editorial board links. We opted for manual resolution of URLs when necessary because it was not onerous, though this process could have been automatized with further use of MTurk HITs.

We obtained editorial board members by launching another batch of HITs. For these HITs, we provided Turkers with the journal title and editorial board URL. We asked Turkers to provide the following information about each editorial board member: first and middle names, last name, title on editorial board (\emph{e.g.}, Editor-in-Chief, Associate Editor, etc.), institutional affiliation, city, state/province, postal code, and country. In the case of first names and middle names, Turkers could provide initials if full names were not available. Turkers could also leave blank other fields if information was not available. Each HIT was initially completed three times. For 229 journals, at least two Turkers reported the same number of editorial board members to within 10\%, and we took these as valid. For these 229 journals, we saved the intersection of editorial board members contained in the valid responses. For the remaining 374 journals, we solicited a fourth response and accepted data as valid based on the same 10\% criterion, again keeping the intersection of board members for each journal. Overall, we obtained data for 435 journals comprising 13067 editorships.

Throughout this paper, we use the word \emph{editorships} rather than \emph{editors} because we take as our fundamental data element a position on an editorial board, and not a unique individual. To restate, we have not de-duplicated the data. If a particular individual serves on the editorial boards of more than one journal, that individual will appear multiple times in our data set, potentially with slightly different representations of their name. For instance, our data set includes ``SS Dragomir'' for \emph{Journal of Inequalities and Applications}, ``Sever S Dragomir'' for \emph{Banach Journal of Mathematical Analysis}, and ``Sever Dragomir'' for \emph{Filomat}. Additionally, as we will discuss further, we intend the term \emph{editorship} to encompass many varying titles on editorial boards.

\subsection*{Coding and cleaning of data}

We coded four variables in our raw data set, used later in our analysis: subfield of journal, publisher, title on editorial board, and country.

As previously mentioned, JCR indexes four categories of mathematical sciences journals, with journals potentially belonging to more than one category. We collapse these four categories into three non-overlapping categories: \emph{pure}, \emph{applied}, and \emph{both}. \emph{Pure} contains journals belonging solely to JCR's \emph{mathematics} category, for example, \emph{Journal of Number Theory}, \emph{Journal of Convex Analysis}, and \emph{Algebra and Logic}.  \emph{Applied} contains journals belonging to one or more of JCR's \emph{mathematics applied}, \emph{mathematical and computational biology}, and \emph{mathematics interdisciplinary applications} categories, but not to \emph{mathematics}. Examples include \emph{Journal of Mathematical Psychology}, \emph{Chaos}, and \emph{Computational Mechanics}. Finally, \emph{both} contains journals that belong to \emph{mathematics} and at least one of JCR's remaining three categories. Examples include \emph{Advances in Calculus of Variations}, \emph{Topology and its Applications}, and \emph{Journal of Pure and Applied Algebra}.

To identify the publisher associated with each editorship, we scraped from the Web a database of journal titles providing each one's ISSN number and publisher \cite{ISSNlist} and joined it to our database of editorships, obtaining a publisher for 100\% of those editorships. Then, we manually aggregated publisher names known to represent the same publisher or be owned by the same parent company. For example, we mapped \emph{Elsevier GMBH}, \emph{Elsevier Sci Ltd}, \emph{Elsevier Science BV}, \emph{Elsevier Science Inc}, and \emph{Elsevier Science SA} all to \emph{Elsevier}, and mapped \emph{Biomed Central Ltd} to \emph{Springer}. This process reduced the total number of publishers in our data set from 156 to 123.

We collapsed the title on editorial board for each editorship into one of three levels: \emph{managing}, \emph{editor}, and \emph{other}. \emph{Managing} captures editorships whose raw titles imply a leadership element, such as editor-in-chief, managing editor, associate managing editor, managing board member, and chief editor. \emph{Editor} comprises raw titles such as editor, associate editor, member of editorial board, editorial committee, and so forth. \emph{Other} comprises raw titles such as honorary editor, advisory board, academic editor, founding board member, and editor emeritus.

To identify the editorship's country of residence, we first used the Google Maps API \cite{Googlemaps} to geocode any country, city, or state Turkers extracted about the editorship. If no geographic information was available or if the geocoding failed, we also tried to geocode the editorship's institutional affiliation. In total, we identified countries for 91.1\% of the editorships.

\subsection*{Initial gender inference}

Gender is a complicated construct touching upon biological sex, social roles, individual gender identity, and more \cite{Lin2015}. Though the most commonly recognized genders are man and woman, other gender identities include transgender, agender, intersex, and many more. The tools we describe below, used to infer gender, are admittedly limiting in that they adopt a binary classification. Still, with this limitation recognized, and with the additional recognition that an individual's gender is most appropriately expressed and explained by that individual, we proceed with gender classification because data is necessary in order to begin quantifying the representation of women. We used two methods to infer gender of editorial board members. 

First, we used \texttt{genderize.io}. This product uses a large corpus of first names and known genders gathered from social networks in order to predict gender. At the time we ran our study, the corpus contained 216286 distinct names across 79 countries and 89 languages. For each first name provided, \texttt{genderize.io} returns a count of the number of times that name appears in the corpus, and corresponding probabilities of (binary) gender based on frequency counts. Intuitively, one has more confidence in the predictions of names with high counts in the corpus, and less confidence in low-count names. To weigh the confidence of high-count names more, we created a modified probability score based on an \emph{a priori} estimate that any given name has equal probability of being a man's or woman's. More specifically, using the probability $p$ and count $c$ for the most likely gender output by \texttt{genderize.io}, we calculated a modified probability,
\begin{equation*}
p_{mod} = \frac{pc+2}{c+4}.
\end{equation*}
We chose the one free parameter above, namely $2$, via a modest amount of manual experimentation. We then accepted the predicted gender for any editorship having $p_{mod} > 0.85$. To verify this as an acceptable threshold, we manually checked names near the margin. Example names having $0.85 \leq p_{mod} \leq 0.855$ include Aman, Dorian, Gerry, Raj, and Tristan, names for which we agree with the most likely predicted gender. In all, our criterion accepted 6555 out of 13067 gender predictions from \texttt{genderize.io}, accounting for 50.2\% of our data set. Within this group, the average inferrence score was $p_{mod} = 0.973$ for women, $p_{mod} = 0.976$ for men, and $p_{mod} = 0.975$ overall.

For the remaining 6512 entries, accounting for 49.8\% of our data set, we inferred gender using MTurk. Each entry was coded by five different Turkers. Each Turker independently located the editor using a search engine based on the information we had previously collected (first and last name, institution, journal name, and so forth). If possible, the Turker then identified the editor a man, a woman, or non-binary. The Turker could also indicate that they could not determine gender. By use of the term \emph{cannot determine}, we mean not that an individual's gender is actually undetermined, but merely that the Turker could not determine it. Less than $0.1$\% of responses were non-binary, representing 17 distinct individuals for each of whom one out the five assigned Turkers selected a non-binary gender. We assigned these Turk responses to cannot determine, meaning that we could not determine how to appropriately place them on the admittedly limiting binary man/woman scale to which we eventually constrain our data. We also asked Turkers to indicate their confidence in the gender they inferred for the editor on a three point scale (1, 2, or 3).

To reconcile the five gender inferences for editorship into a single inference, we created an aggregate gender score. More specifically, we assigned low, middle, and high confidence inferences of women the values of $+1/3$, $+2/3$, and $+1$ respectively. Inferences of men similarly received values of $-1/3$, $-2/3$, and $-1$, and cannot determine responses received a value of zero. We then averaged the five Turker responses to generate a final gender score that ranged from $+1.0$ (very likely a woman) to $-1.0$ (very likely a man). We experimented with more sophisticated machine learning techniques such as random forests to infer gender using the five Turker responses and confidence levels, but we found the simple gender score approach was equally accurate and more transparent.  Of the 6512 editorships handled via MTurk, 99\% of them had nonzero gender scores, and we inferred a gender depending on whether the score was negative or positive. The remaining 1\% of editorships ($n=67$) had a gender score of zero and are listed as undetermined in our results.

\subsection*{Validation and calibration of inferred gender}

We performed an additional calibration and validation step to estimate the accuracy of our gender inference process, and critically, to ensure the process had no distributional biases that estimated too many men or women overall. To do so, we collected expert-inferred genders for a stratified sample of editors, as we now explain.

We include one sample stratum for editorships inferred as women by \texttt{genderize.io}, one sample stratum for editorships inferred as men by \texttt{genderize.io}, and one sample stratum for each interval of MTurk aggregate gender score in the sets -1, (-1,-0.8], (-0.8,0.6], (-0.6,-0.4], (-0.4,-0.2], (-0.2,0.2), [0.2,0.4), [0.4,0.6), [0.6,0.8), [0.8,1), and 1. To summarize, we had 13 strata: two sampled from \texttt{genderize.io} and 11 from across MTurk gender scores. For each of the 13 strata, we randomly sampled 30 editorships and ourselves served as the experts. One of us first inferred the gender of the editorship using a search engine. Validation samples were shuffled so that we did not know which stratum an editorship belonged to, and we selected \emph{cannot be determined} if we could not robustly infer an editorship's gender. When the expert-inferred gender agreed with the gender as originally inferred by \texttt{genderize.io} or MTurk, we took the result as correct. When the two disagreed, the remaining expert resolved the disagreement.

\begin{figure}[!h]
\includegraphics[width=\textwidth]{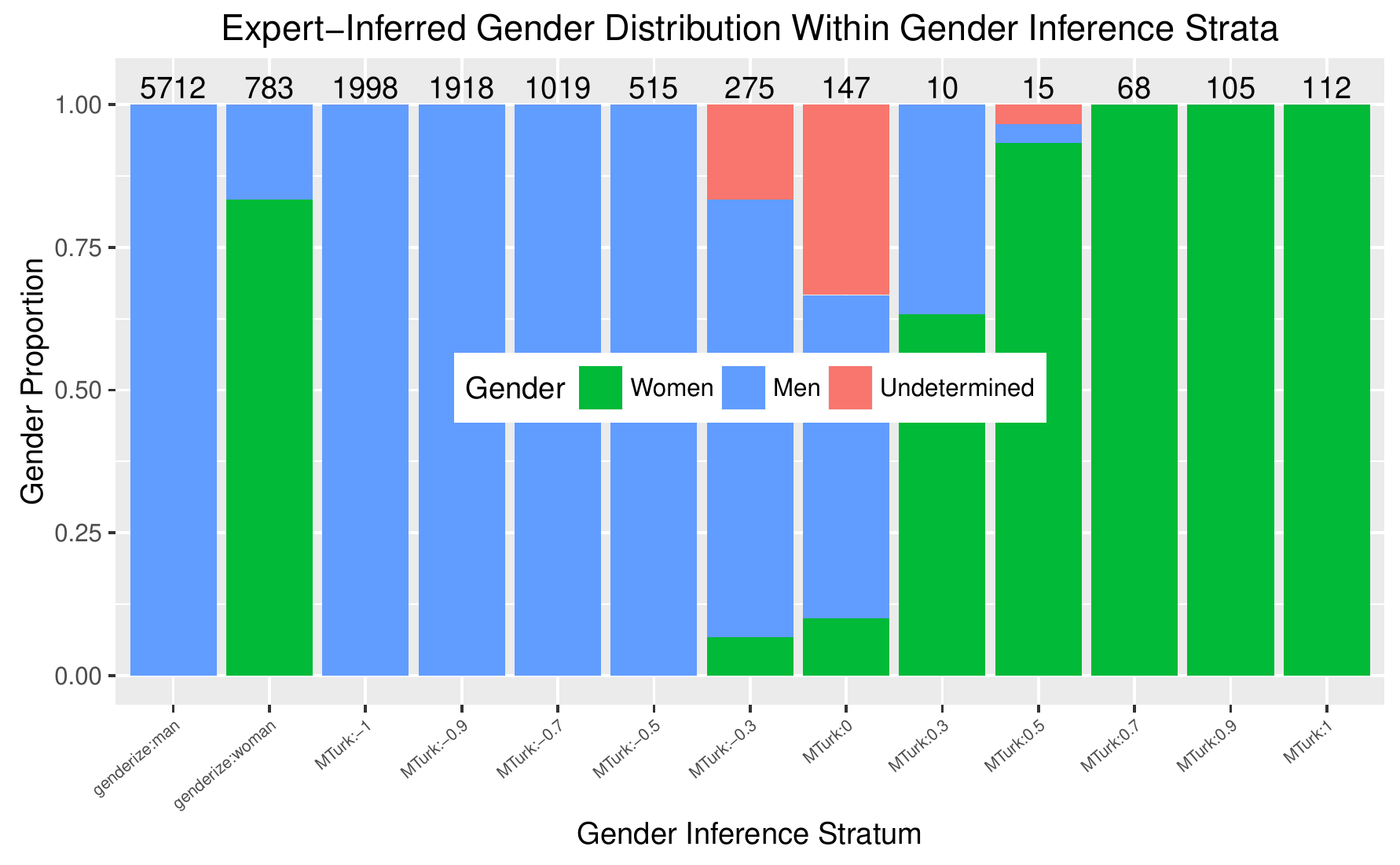}
\caption{{\bf Distribution of expert-inferred genders for each stratum in gender calibration procedure}.  Each bar shows the distribution of 30 expert gender inferences performed by at least one of the investigators. The 30 editorships comprising each bar are sampled randomly from the appropriate stratum of our data set. The two left bars show the distribution of expert-inferred genders for editorships processed through the \texttt{genderize.io} gender inference tool. The eleven rightmost bars show the distribution of expert-inferred genders for different levels of gender score compiled from our Amazon Mechanical Turk (MTurk) data. The number in each stratum label refers to the midpoint of the gender score interval. For instance, \emph{MTurk:-1} corresponds to editorships inferred by all five Turkers to be men, while \emph{MTurk:0.7} corresponds to editorships with MTurk gender scores in the interval $[0.6,0.8)$, likely women; see text for further explanation. The number on top of each bar gives the number of unvalidated editorships, that is, the number of editorships in our full data set falling within the stratum, excluding the 30 that were randomly sampled.}
\label{fig3}
\end{figure}

Fig~\ref{fig3} shows the results of the validation procedure. Each bar visualizes the distribution of expert-inferred genders for the 30 samples within the stratum. The two left-most bars show the results for names inferred to be women and men by \texttt{genderize.io}. The remaining 11 strata show results for the aforementioned intervals of MTurk gender score. First, expert-inferred genders agreed with the genderize.io result for 25 of the 30 (83.3\%) editorships in the woman \texttt{genderize.io} strata, and all 30 names inferred as men agreed. The higher percentage of incorrect women's name inferences may reflect the mismatch between the gender distribution that the \texttt{genderize.io} algorithm was trained on and the gender distribution most often seen among mathematical sciences editorships. The algorithm is trained on social media users, whose genders are close to evenly distributed, while the mathematical sciences editorships to which we applied the algorithm are overwhelmingly men. Second, the gender scores that result from Turker responses correlate with expert-inferred genders. Scores $\le -0.4$ and $\ge +0.6$ exhibited no disagreements. We did find two errors (6.7\%) in scores within the $+0.4$ stratum. While more discrepancies occurred at the $+0.2$, $0.0$, and $-0.2$ levels, these strata are rare within our dataset (see numbers on top of bars in figure) accounting for a combined 4.0\% of all editorships. Many of the first names in this range of MTurk gender scores are solely initials, \emph{e.g.}, ``L'' and ``AJ,'' while other examples include ``Jan,'' ``Lee,'' and ``Khvicha.'' 

Finally, we attempted to eliminate structural bias in our results by computing well-calibrated \cite{NicCar2005} gender probabilities based on the analysis of sampled strata, above. For each editorship, we used the strata validations to introduce calibrated $[P_w, P_m, P_u]$ expected gender distributions that sum to $1.0$, where $P_w$ is the probability of being a woman, $P_m$ is the probability of being a man, and $P_u$ is the probability of being other or undetermined. For the 390 editorships in our sampled strata, we assigned distributions of $[1,0,0]$, $[0,1,0]$, or $[0,0,1]$ because we believed the expert-driven procedure eliminated gender inference errors. For each of the remaining 12667 unvalidated editorships, we assigned the editorship the observed distribution of expert-inferred gender probabilities within its strata. For example, all 783 unvalidated editorships that were inferred to be women's names by \texttt{genderize.io} received a corrected gender distribution of $[0.833, 0.167, 0.0]$ corresponding to the second bar from the left in Fig~\ref{fig3}.

We used our inferred gender probabilities to estimate the accuracy of our raw (uncalibrated) inference strategy. For example, in the stratum inferred as women by \texttt{genderize.io}, the accuracy of the raw inferences can be estimated at $83.3\%$ (25/30). By summing our calibrated gender probabilities across the full dataset, we can compute an expected confusion matrix, shown in Table~\ref{table1}. The matrix estimates that overall, our raw gender inference process is correct 97.5\% of the time (sum of diagonal elements). Crucially, we expect the actual number of women to be approximately 1\% lower than the raw estimates (8.9\% vs 9.9\%), primarily due to the overestimate of women by the \texttt{genderize.io} API. In addition, despite the relatively small number of women in the dataset, we can estimate that 97\% of the actual women are inferred to be women, so that the sensitivity (or recall) for women is $8.6/8.9 = 97\%$. Conversely, we estimate that 89\% of editorships inferred to be women are actually women, so that the specificity (or precision) for women is $8.6/9.9=89\%$. 

\begin{table}[!ht]
\centering
\caption{
{\bf Estimated confusion matrix comparing gender inferences made by our raw procedure (\texttt{genderize.io} and MTurk) with those made by our calibrated inference strategy, which is based on expert analysis of sampled strata. This matrix estimates that our raw procedure is correct 97.5\% of the time (sum of the diagonal elements). However, we expect the actual number of women to be approximately 1\% lower than the raw estimates (8.9\% vs 9.9\%), primarily due to the overestimate of women by the \texttt{genderize.io} API. Throughout the remainder of our study, we use the calibrated gender inferences.}}
\begin{tabular}{|r|c|c|c|c|}
\hline
& Expert woman & Expert man & Expert other & Total\\ \hline
Inferred woman & 8.6\%& 1.3\%& 0.1\% & 9.9\%\\ \hline
Inferred man & 0.2\%& 88.7\% & 0.6\% & 89.6\%\\ \hline
Inferred other & 0.1\%& 0.2\%& 0.2\% & 0.5\% \\ \hline
Total & 8.9\%& 90.3\%& 8.5\% &  \\ \hline
\end{tabular}
\label{table1}
\end{table}

The validation and calibration procedures suggest our gender inference is robust, and we use the calibrated gender probabilities presented above throughout the remainder of this paper. As we note in the following section, our bias correction and accuracy estimates are limited by the accuracy and representativeness of the expert-inferred genders. While we believe these to be correct and note that the sample size is large enough to rule out serious issues, some miscalibrations may remain. 

\subsection*{Limitations}

We provide several points of caution for the reader. First, because we began with 605 journals from JCR and accepted only those for which at least two Turkers found sufficiently close lists of editorships, our data set must be viewed as a convenience sample of the overall landscape. 

Second, as previously discussed, the gender proportions reported below are calibrated values that result from an inference process. Therefore, gender proportions for some journals may imply non-integer editor counts for each gender (though the totals always sum to integers).  We must be especially cautious for names inferred to be women by \texttt{generize.io} both because of the error we observed within the validation sample and because of the relatively large size of this stratum. Our adjustment methodology may misreport representation of women on journals that have many editorships inferred to be women by \texttt{genderize.io}. For example, \textit{Envirometrics}, one of the journals with highest representation of women, had a raw estimate of 33.3\% women. The calibration procedure lowered this number to 28.4\% because it predicted some of the 19 editorships that were inferred to be women based on first name were actually men. However, after communicating with the editors of \textit{Envirometrics} we learned the raw inferences were indeed more accurate, and the actual percentage of women is 34.5\%. Despite these caveats, across all journals, the mean difference between the raw estimate of the percentage of women and the calibrated estimate is 1\%, and overall, we expect calibrated values to be more accurate, as explained in the previously. Our publicly available data set (see Data Availability Statement) contains both raw and calibrated inferences. 

Third, some miscalibrations no doubt remain in the data. An individual expert or team of experts with sufficient knowledge might be able to outperform our gender inference process for any particular journal editorial board. However, part of the purpose of our study is to propose a methodology sufficiently automatic that it can be used at a very large scale. 

Fourth, and finally, these data were gathered during early 2016, and are a snapshot of that time. Online editorial board information may change at any time, even between data gathering and the publication of this manuscript.

\section*{Results}
\label{sec:resultsanddiscussion}

\subsection*{Overview}

Our data set of 13067 editorships is 8.9\% women, 90.3\% men, and 0.8\% undetermined. The coded titles of these editorships are 82.2\% \emph{editor}, 5.1\% \emph{managing}, and 12.7\% \emph{other}.

There are 435 journals belonging to 123 publishers. The sizes of the editorial boards form a distribution that is very roughly bell-shaped (not shown) with a mean of 30 and a median of 29. Of the journals, 35.6\% are \emph{pure}, 43.9\% are \emph{applied}, and 20.5\% are \emph{both}. If we instead consider the number of editorships belonging to each subfield, we find 27.7\% are \emph{pure}, 51.9\% are \emph{applied}, and 20.4\% are \emph{both}. The distribution of five-year impact factor is strongly right-tailed, with a median of 0.84 and a mean of 1.1.

For 91.1\% of the data, we have a country of residence for the editorship, with 86 countries appearing. The ten countries contributing the greatest number of editors are the US (33.6\%), Great Britain (7.4\%), France (6.7\%), Germany (6.6\%), Italy (4.5\%), Canada (4.0\%), Japan (3.8\%), China (3.8\%), Russia (2.4\%), and Australia (2.0\%), where percentages are with reference to the large subset for which country data is available.

In the following subsections, we examine gender representation by journal, journal subfield, impact factor, publisher, editorship title, and editorship country. Throughout, we take as our primary goal the quantification of gender representation and the identification of outliers, that is, any subgroups that have statistically significantly low or high representation of women. In the Conclusions section, we pose a number of questions raised by our results.

\subsection*{Gender analysis by journal}

\begin{figure}[!h]
\includegraphics[width=\textwidth]{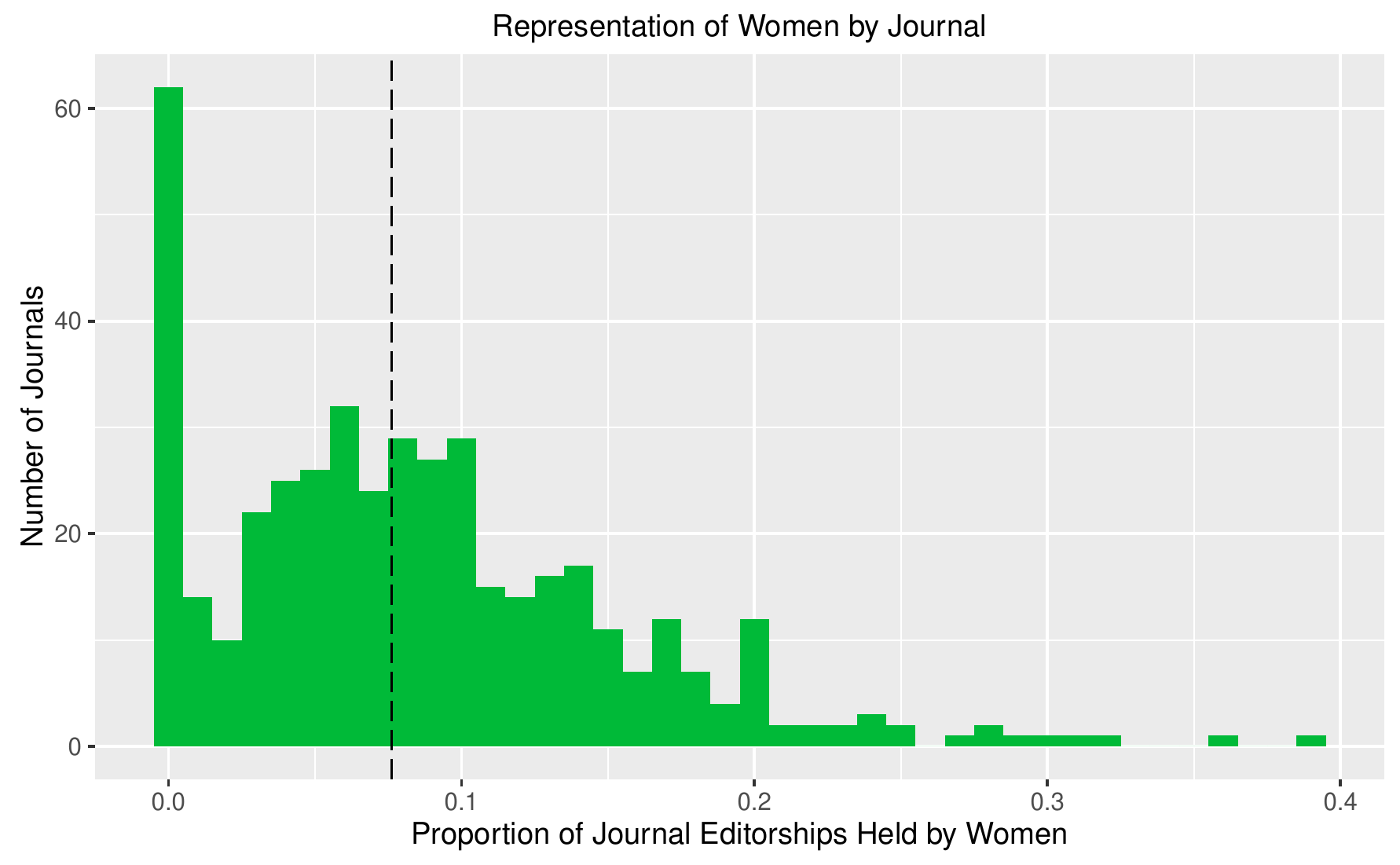}
\caption{{\bf Distribution of journals by proportion of editorships held by women}.  The median journal is 7.6\% women, indicated by the vertical dashed line. The large spike at the left of the histogram represents 62 journals having less than one-half percent women on the editorial board. In fact, 51 journals have no women, accounting for 11.8\% of the 435 journals in our data set.}
\label{fig4}
\end{figure}

Fig~\ref{fig4} shows the distribution of journals by proportion of editorships held by women. Overall, the median journal has 7.6\% women editorships as indicated by the vertical dashed line. The large spike at the left of the histogram represents 62 journals having less than one-half percent women. In fact, 51 journals have no women, accounting for 11.7\% of the journals in our data set. A few examples include \emph{Annals of Mathematics}, \emph{Communications on Pure and Applied Mathematics}, \emph{Inventiones Mathematicae}, \emph{Journal of Algebraic Geometry}, \emph{Journal of Differential Geometry}, and \emph{Mathematische Zeitschrift}. Of these 51 journals, 29 are \emph{pure}, 11 are \emph{applied}, and 11 are \emph{both}. We return to measure variation in mathematical areas in the following section.

\begin{figure}[!h]
\includegraphics[width=\textwidth]{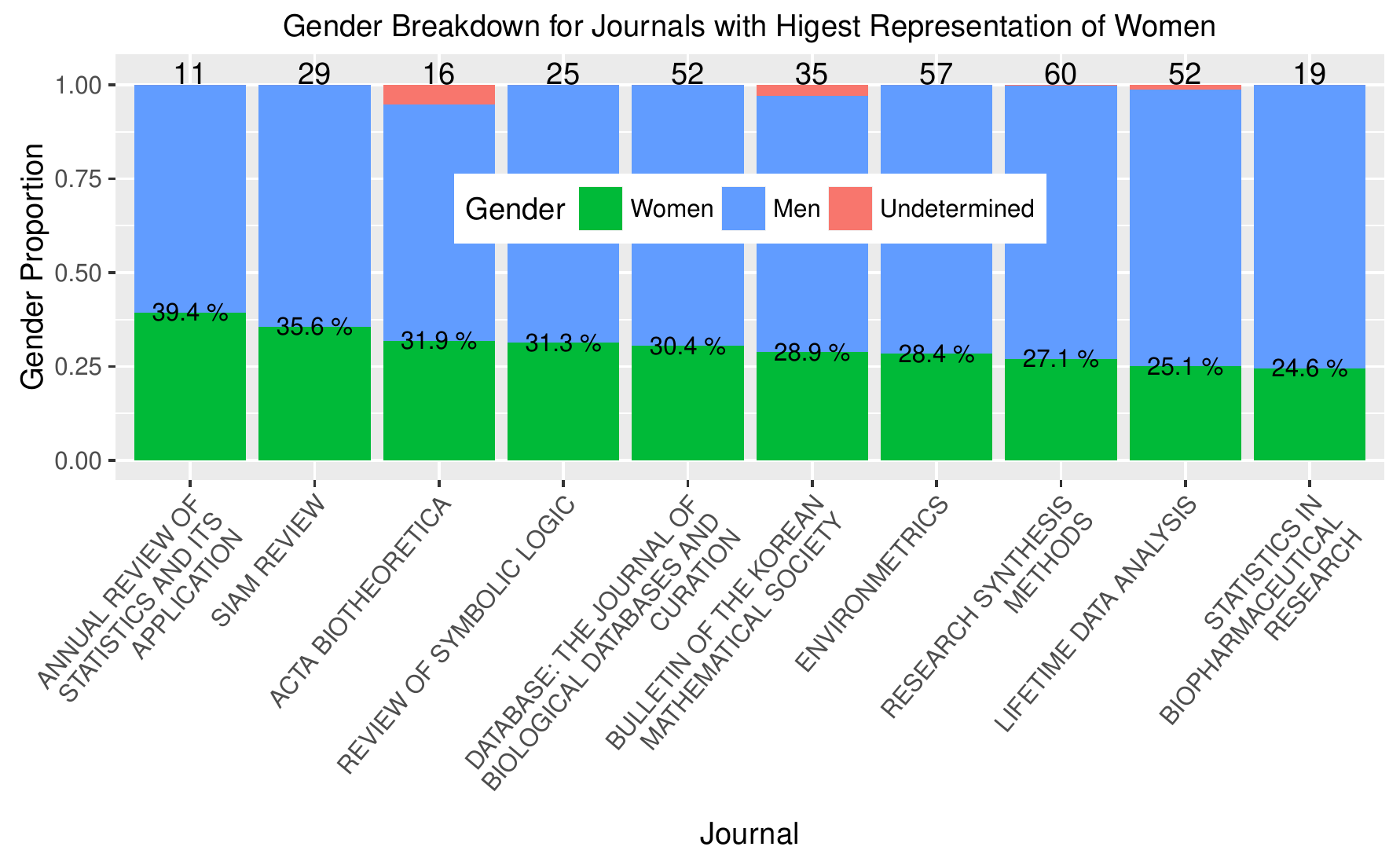}
\caption{{\bf Journals with highest representation of women.} The number on top of each bar gives the quantity of editorships on each journal and the percentages within the bars state the proportion held by women, represented by the green area. Of the top 10 journals, only \emph{Bulletin of the Korean Mathematical Society} is in the subfield category \emph{pure}. In this figure, we exclude editorial boards with fewer than ten editors; there is only one such journal that would have been in the top ten, namely \emph{Periodica Mathematica Hungarica}.}
\label{fig5}
\end{figure}

Fig~\ref{fig5} shows the gender breakdown of the ten editorial boards with highest representation of women. In this analysis, we exclude journals with fewer than ten editorships. (There is only one such journal that would have been in the top ten, namely \emph{Periodica Mathematica Hungarica}, which has only three editorships). The three journals with the highest representation are \emph{Annual Review of Statistics and its Applications} (39.4\%), \emph{SIAM Review} (35.6\%), and Acta Biotheoretica (31.9\%). Of the top ten journals, only \emph{Bulletin of the Korean Mathematical Society} is in the subfield category \emph{pure}; the remaining nine are either \emph{applied} or \emph{both}.

For each journal, we perform a chi-squared ($\chi^2$) test for difference in proportions to compare the percentage of women on that journal to the percentage of women in the remaining journals grouped together. This test determines whether that journal has representation of women that is statistically significantly different from the background level. Whether a particular journal turns out to be statistically significant depends not just on the percentage of women, but on the size of the whole editorial board; for instance, a lower proportion of women is more significant for a large editorial board than for a small one. For all $\chi^2$ tests through this manuscript, we report $p$-values at or below the $0.05$ significance threshold. In the present subsection, we use Holm-adjusted $p$-values which account for the large number of hypotheses tested (each of the 435 journals against the remaining background) \cite{Hol1979}. From left to right in Fig~\ref{fig5}, the five journals that are significantly different from the background level are \emph{SIAM Review} ($p = 0.003$), \emph{Database} ($p < 0.001$), \emph{Environmetrics} ($p < 0.001$), \emph{Research Synthesis Methods} ($p = 0.002$), and \emph{Lifetime Data Analysis} ($p = 0.05$). All five journals are classified as \emph{applied}. The $\chi^2$ test reveals no journals with significantly lower representation of women than the background level.

\subsection*{Gender analysis by journal subfield}

Within the pool of editorships on journals classified as \emph{applied}, representation of women is 10.3\%, while the percentages are 7.2\% for \emph{pure} and 7.4\% for \emph{both}. For our analysis of journals above, we performed $\chi^2$ tests to compare each journal to the background level of the remaining journals grouped together, and later, we will do this for publishers. Presently, in the case of journal subfield, there are only three groups, so we instead perform three pairwise $\chi^2$ tests. We find that \emph{applied} has higher representation than \emph{pure} ($p<0.001$) and \emph{applied} has higher representation than \emph{both} ($p<0.001$). The $\chi^2$ test could not detect a significant difference in representation of women between \emph{pure} and \emph{both}. Therefore, it appears that journals with a pure component are associated with lower representation of women than those that are exclusively applied.

\subsection*{Gender analysis by impact factor}

For the 435 journals in our data set, 416 of them have available a five-year impact factor. To examine possible associations between a journal's impact factor and the percentage of its editorships held by women, we produce a scatter plot of these two variables. This plot (not shown) reveals no obvious trend. To assess a possible association more quantitatively, we compute Spearman's rank correlation coefficient $r_s$, which measures the degree to which the relationship two variables is monotonic, without assuming linearity. In our case, we find a significant correlation of $r_s = 0.14$ ($p=0.004$).

Within our data set, the journals published by SIAM Publications have amongst the highest representation of women, and additionally, have high impact factors. As an additional test, to assess whether SIAM journals are responsible for the degree of correlation that we found above, we excluded the seven SIAM journals in our data set and recalculated $r_s$ for the remaining 409 journals. The associated $r_s = 0.11$, is still significant ($p=0.02$). Overall, we find a statistically significant but weak positive association between a journal's impact factor and the percentage of its editorships held by women.

\subsection*{Gender analysis by publisher}

\begin{figure}[!h]
\includegraphics[width=\textwidth]{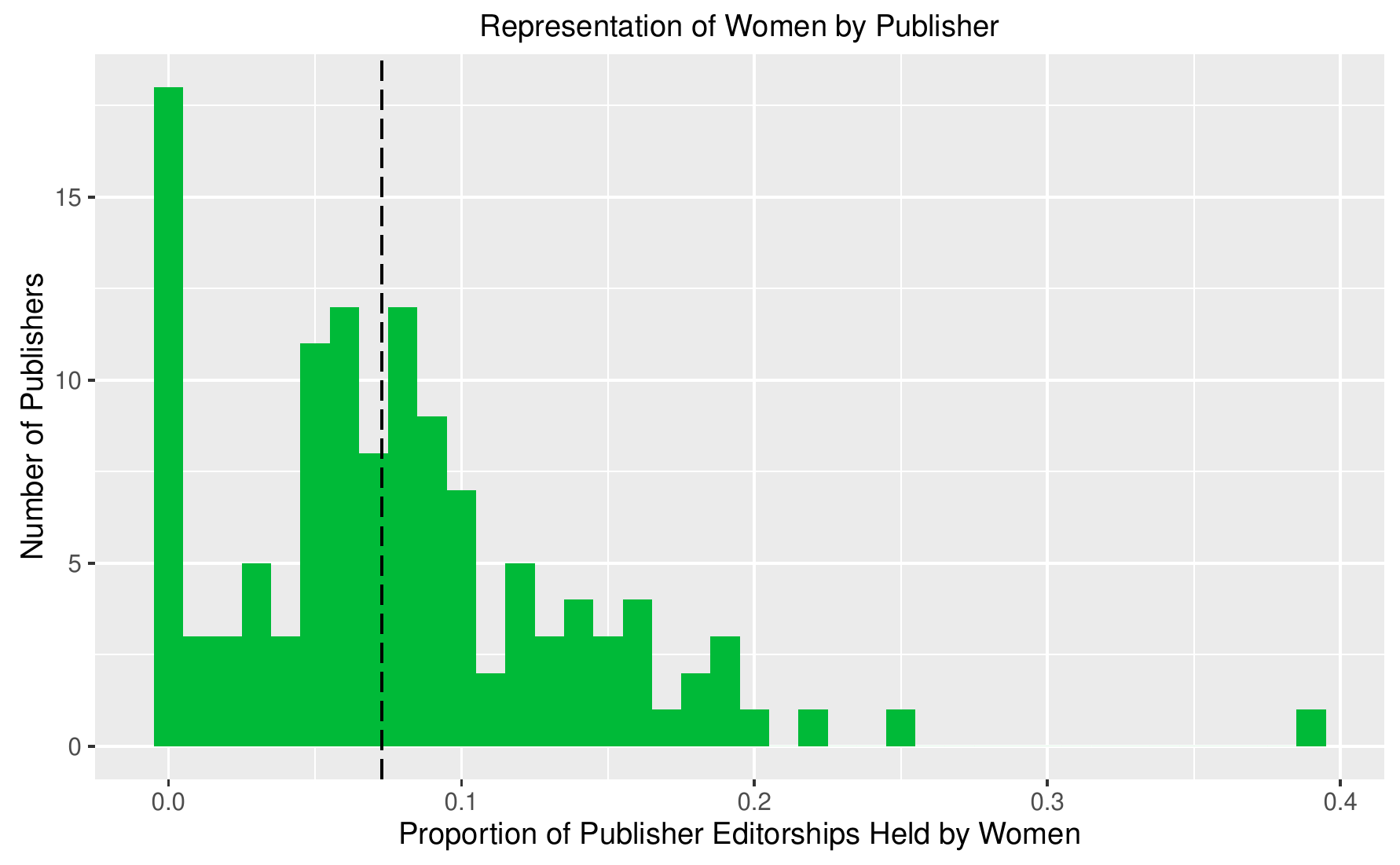}
\caption{{\bf Distribution of publishers by proportion of editorships held by women}.  The median publisher has 7.3\% editorships held by women, indicated by the vertical dashed line. The large spike at the left of the histogram represents 18 publishers having less than one-half percent women on the editorial board. Of these, 16 publishers have no women at all. All 18 are all publishers with only one journal appearing in our data set. The three, small rightmost bars also correspond to publishers each with a single journal each in our data set, the largest of which has 21 editorships.}
\label{fig6}
\end{figure}

Fig~\ref{fig6} shows the distribution of publishers by proportion of editorships held by women. We pool together all journals in our sample belonging to each publisher. Overall, the median publisher has 7.3\% editorships belonging to women, as indicated by the vertical dashed line. The large spike at the left of the histogram represents 18 publishers having less than one-half percent editorships belonging to women. Of these, 16 publishers have no women editorships at all. All 18 are all publishers with only one journal appearing in our data set, and the largest such journal/publisher has 39 editorships. At the other end of the histogram, The three, small rightmost bars in the histogram correspond to (from right to left) \emph{Annual Reviews}, \emph{American Statistical Association}, and \emph{Ramanujan Mathematical Society}. These three publishers also have a single journal each in our data set, the largest of which has 21 editorships. 

\begin{figure}[!h]
\includegraphics[width=\textwidth]{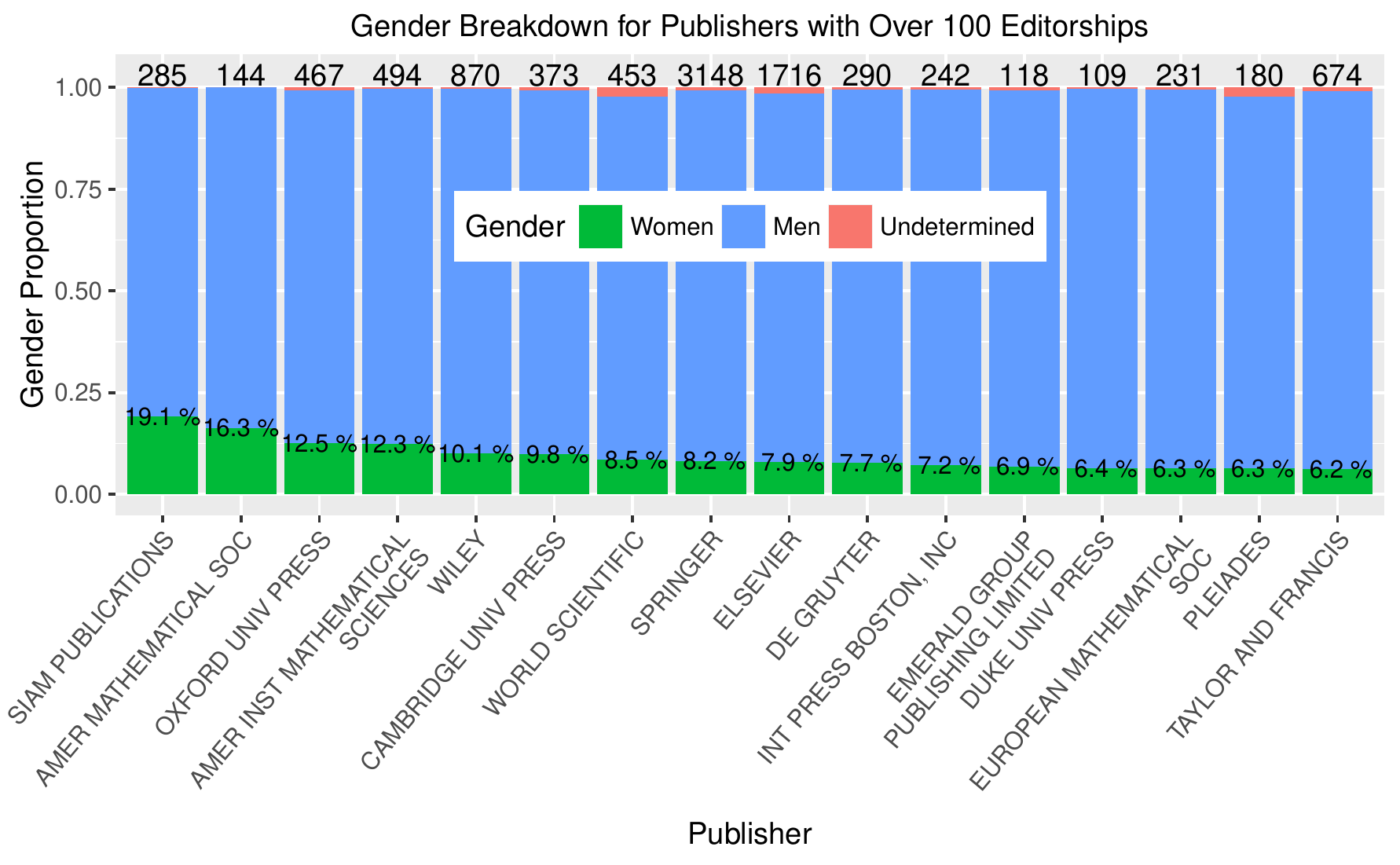}
\caption{{\bf Publishers with over 100 editorships, ordered by representation of women.} The number on top of each bar gives the quantity of editorships from each publisher and the percentages within the bars state the proportion held by women, represented by the green area. The four publishers with the highest representation of women are scientific societies or university publishing houses.}
\label{fig7}
\end{figure}

We now focus on larger publishers, which we define as publishers having at least 100 editorships in our dataset. There are 16 such publishers; Fig~\ref{fig7} shows their gender breakdown. The four high-editorship publishers with the greatest representation of women are \emph{SIAM Publications} (19.9\%), \emph{American Mathematical Society} (16.3\%), \emph{Oxford University Press} (12.5\%), and \emph{American Institute of Mathematical Sciences} (12.3\%). These four publishers are either scientific societies or university publishing houses. The commercial publishing houses with the highest representation of women are Wiley (10.0\%), World Scientific (8.5\%), Springer (8.2\%), and Elsevier (7.9\%). The last three are close to the median value in the histogram of Fig~\ref{fig6} (within one or two histogram bins).

For each publisher, we perform a $\chi^2$ test for difference in proportions to compare the percentage of women within that publisher to the percentage women within all other publishers grouped together, again using Holm-adjusted $p$-values to account for multiple comparisons. Whether a particular publisher turns out to be statistically significant depends not just on the percentage of women, but on the size of the pool of editorships within that publisher. The only publisher that is significantly higher than the background level is SIAM Publications ($p < 0.001$). The $\chi^2$ test reveals no publishers with significantly lower representation than the background level.

\subsection*{Gender analysis by editorship title}

Within the pool of editorships with editorial board titles classified as \emph{editor}, $9.0\%$ are women, while the percentages are 7.4\% for \emph{managing} and 8.7\% for \emph{other}. As with our gender analysis by journal subfield, since there are only three groups of editorship titles, we perform pairwise $\chi^2$ tests. These tests do not detect any differences in representation of women. However, a lower representation of women at higher levels of leadership is consistent with previous research (for example, \cite{Stewart2007}). Because a great variety of editorial board titles are used, and because our coding of this variable into three levels is admittedly coarse, a more nuanced coding might reveal deeper insights.

\subsection*{Gender analysis by country}

\begin{figure}[!h]
\includegraphics[width=\textwidth]{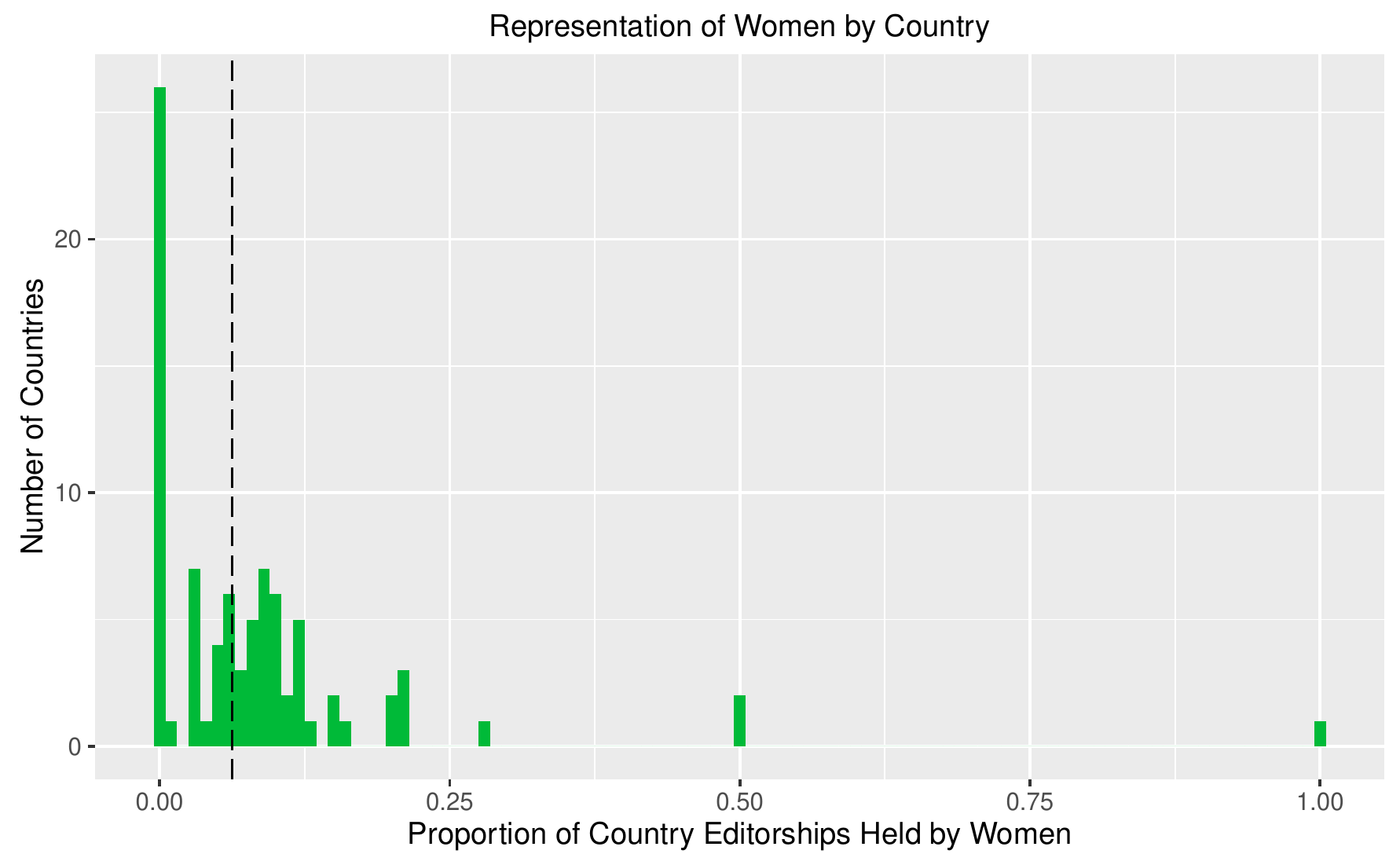}
\caption{{\bf Distribution of countries by proportion of editorships held by women}.  The median country has 6.3\% editorships held by women, indicated by the vertical dashed line. The large spike at the left of the histogram represents 26 countries with no women in our data set. However, all 26 of these countries have at most 3 editorships per country. At the other end of the histogram, the four rightmost bars, representing 6 countries, each have at most 5 editorships.}
\label{fig8}
\end{figure}

\begin{figure}[!h]
\includegraphics[width=\textwidth]{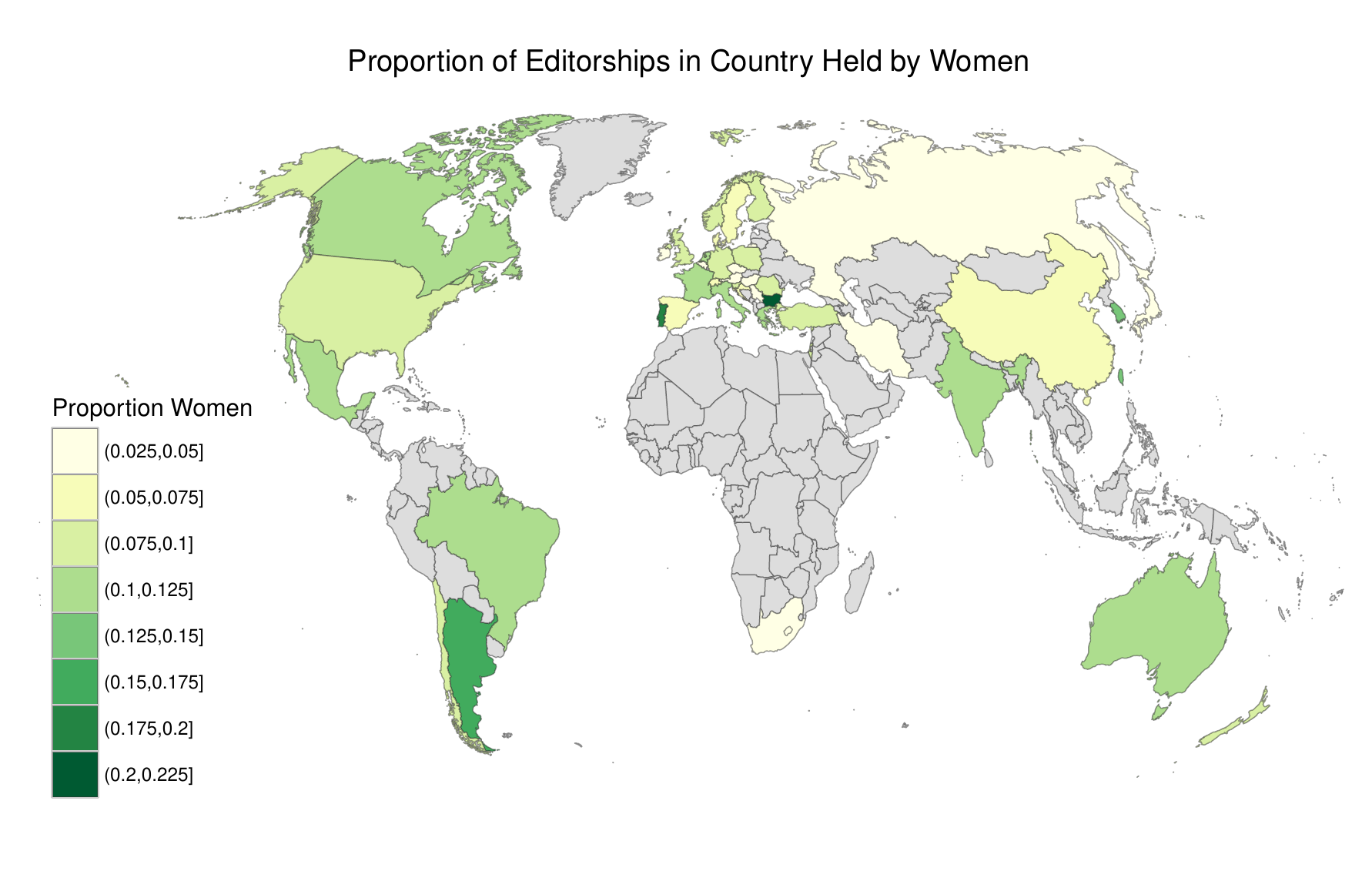}
\caption{{\bf Proportion of editorships in each country that are held by women.} Countries with no editors are shown in gray.}
\label{fig9}
\end{figure}

Fig~\ref{fig8} shows the distribution of countries by proportion of editorships held by women. Overall, the median country has 6.3\% of editorships held by women, as indicated by the vertical dashed line. The large spike at the left of the histogram represents 26 countries with no women in our pool of editorships. However, all 26 of these countries have at most 3 editorships per country. At the other end of the histogram, the four rightmost bars (representing 6 countries) each have at most 5 editorships. Fig~\ref{fig9} visualizes our data on a world map, with the percentage of editorships held by women indicated by shading.

\begin{figure}[!h]
\includegraphics[width=\textwidth]{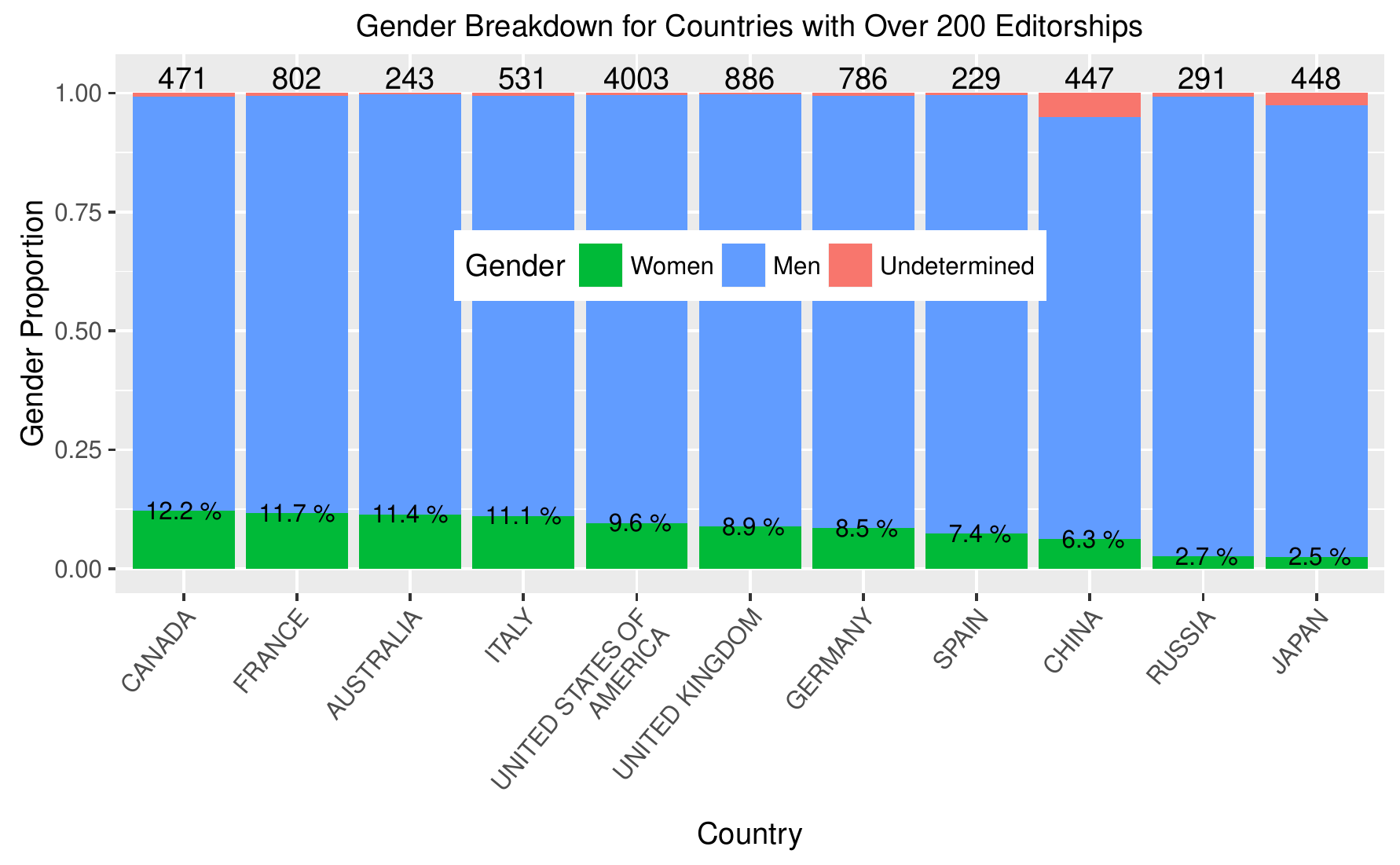}
\caption{{\bf Countries with over 200 editorships, ordered by representation of women.} The number on top of each bar gives the quantity of editorships from each country and the percentages within the bars state the proportion held by women, represented by the green area.}
\label{fig10}
\end{figure}

We now focus on countries comprising larger numbers of editorships; we choose a threshold of 200. There are 11 such countries and Fig~\ref{fig10} shows their gender breakdown. The four high-editorship countries with the greatest representation of women are Canada (12.2\%), France (11.7\%), Australia (11.4\%), and Italy (11.1\%). The United States (9.6\%) is fifth within the group of high-editorship countries, and 24th overall.

For each country, we perform a $\chi^2$ test for difference in proportions to compare the percentage of editorships in residing in that country that are held by women to the percentage in the remaining countries grouped together, again using Holm-adjusted $p$-values to account for multiple comparisons. Whether a particular country turns out to be statistically significant depends not just on the percentage of women, but on the size of the pool of editorships within that country. The only countries that are significantly different from the background level are Russia ($p = 0.03$) and Japan ($p < 0.001$), both with lower representation.

\section*{Conclusion}
\label{sec:conclusion}

Within the United States, women account for approximately 51\% of the population, 42\% of bachelor's degrees in the mathematical sciences, 29\% of doctoral degrees in the mathematical sciences, and 15\% of tenure-stream faculty at doctoral granting mathematical sciences departments \cite{Census2016,IPEDSCompletions2016,IPEDSSalaries2016,SED2016,AMS2016}. From our present study, we now know that women account for merely 8.9\% of mathematical sciences journal editorships, with the median journal having 7.6\% editorships held by women. We therefore conclude that the degree of underrepresentation on mathematical sciences journal editorial boards is even more severe than in the field at large.

We have measured and reported on gender representation, and our results raise many questions. For instance, while we find additional underrepresentation at the editorial board level, we have not shed light on the processes responsible. Are women faculty not being sufficiently considered for editorial boards? Are they being considered, but in the end not being asked to serve? Are they being asked to serve but choosing not to? Mathematical sciences journals might benefit from attention to these questions. It might also be useful to disaggregate the role other factors play in editorship selection, including the Carnegie classification of an editor's institution (\emph{e.g.}, R1) and the editor's seniority. Finally, recall that our data set was not de-duplicated, meaning that some individuals appear multiple times, once for each journal they serve on, potentially with slight variations in the presentation of their name. Further work could de-duplicate the data set, which would facilitate an investigation of the network structure of editorships; a previous network study has yielded insight into gender representation in computer science faculty hiring \cite{WayLarCla2016}. De-duplication would also let us assess the extent to which individual women serve in multiple editorial roles.

In the Introduction, we mentioned that serving on an editorial board can contribute positively to professional advancement and networking, and that women excluded from editorial board membership do not have access to these advantages. Some have argued to us, however, that one or two women serving on an editorial board comprised largely of men -- despite these women's actual strong qualifications -- might be assumed by onlookers to have obtained their positions precisely because of a journal's desire to diversify its editorial board. If this is indeed the public's perception, then women editors might not be receiving the same advantages of professional status. These assertions would benefit from formal study.

Despite our grim findings for the field at large, certain groups stand out as having higher representation of women than the background level. These include five journals in the applied mathematical sciences and a publishing house within an applied mathematics professional society. The background statistics we report on women's representation at the bachelor's, doctoral, and, faculty levels are for fields in the mathematical sciences aggregated together. Further work might make careful comparisons between the editorial board gender measurements and women's representation in specific mathematical subfields. Similarly, for the countries revealed by our statistical tests to have low representation of women, it would be worthwhile to assess the representation in relation to women's representation within each country's academic ranks.

After our study was complete (and only then), we reached out to the aforementioned five journals and one professional society to ask for their perceptions of causes for their higher representation of women. Responses included: a higher representation of women in particular subfields; specific attention to women's representation by leadership; a strong desire to ``just to have the best people possible'' on the editorial board, which includes women; and sensitivity of men to issues of sexism and underrepresentation. While our present work cannot confirm or deny these explanations, we nonetheless hope that we contribute to raising awareness of gender underrepresentation, and possible ways to address it.

In carrying out our work, we have presented a methodology that is semi-automated and scalable. We recognize the limitations of the man/woman gender scale we have used, and stress the importance of developing better ways of measuring and discussing gender (or lack thereof) that give appropriate accuracy, recognition, and dignity to all individuals and communities. Equally, the complexity of the gender inference task reminds us that the best data on gender comes from individuals themselves, highlighting a benefit that would be obtained if publishers asked editors to self-report gender within an appropriately flexible scale. Still, we hope our techniques might be refined and used in the future to take additional snapshots of the mathematical sciences, shedding light on trends over time. Similarly, our methodology could be applied to fields such as engineering, computer science, and physics, disciplines in which women's representation on editorial boards is also thought to be low.

\section*{Acknowledgments}
This project was inspired by conversations with Linda Chen and Tara Holm. Jude Higdon-Topaz, Sarah Iams, Karen Saxe, and Mary Silber gave us valuable feedback. Andrew Bernoff, Alicia Johnson and Heather Metcalf provided detailed comments and many helpful suggestions on a draft of this manuscript. We are grateful to the dedicated and skilled Amazon Mechanical Turk workers who made possible the efficient completion of this research.

\nolinenumbers



\end{document}